\documentclass[french,a4paper,12pt]{article}
\usepackage{tikz}
\usepackage{epic}
\usepackage[T1]{fontenc}
\usepackage[english]{babel}
\usepackage{amsmath,amsthm}
\usepackage{amssymb}
\usepackage[all]{xy}
\usepackage{stmaryrd}
\usepackage[francais]{layout}
\theoremstyle{definition}
\newtheorem{Thm}{Theorem}

\newtheorem*{coroll}{Corollary}

\newtheorem{lemma}{Lemma}
\theoremstyle{definition}
\newtheorem{Df}{Definition}
\newtheorem{Rmk}{Remark}

\newcommand{\Hom}{\mathop{\mathrm{Hom}}}

\newcommand{\CH}{\mathop{\mathrm{CH}}}
\newcommand{\Ch}{\mathop{\mathrm{Ch}}}

\newcommand{\Corr}{\mathop{\mathrm{Corr}}}

\newcommand{\Spec}{\operatorname{Spec}}
\newcommand{\Mor}{\operatorname{Mor}}

\newcommand{\CM}{\operatorname{CM}}
\newcommand{\CCC}{\operatorname{C}}
\newcommand{\CR}{\operatorname{CR}}
\newcommand{\End}{\operatorname{End}}
\newcommand{\SB}{\operatorname{SB}}
\newcommand{\fonction}[5]{ #1~:\begin{array}{ccc}
 #2 & \longrightarrow & #3 \\
 #4 & \longmapsto & #5 \end{array}}

\author{}
\title{}
\date{}
\begin{document}
\vspace{1cm}
\maketitle

\begin{center}\large{\textbf{A going down theorem for Grothendieck Chow motives}}

\small{Charles De Clercq\footnote{\textbf{Address : }\textsc{Université Paris VI, 4 place Jussieu, 75252 Paris CEDEX 5.}
\textbf{Email :} $\mathsf{\mbox{declercq@math.jussieu.fr}}$}}
\end{center}

\textbf{Abstract.} Let $X$ be a geometrically split, geometrically irreducible variety over a field $F$ satisfying Rost nilpotence principle. Consider a field extension $E/F$ and a finite field $\mathbb{F}$. We provide in this note a motivic tool giving sufficient conditions for so-called outer motives of direct summands of the Chow motive of $X_E$ with coefficients in $\mathbb{F}$ to be lifted to the base field. This going down result has been used S. Garibaldi, V. Petrov and N. Semenov to give a complete classification of the motivic decompositions of projective homogeneous varieties of inner type $E_6$ and to answer a conjecture of Rost and Springer.

\section*{Introduction}

Throughout this note $F$ will be the base field and by an $F$-variety we will mean a smooth, projective scheme over $F$. Given an $F$-variety $X$, we denote by $\Ch(X)$ the Chow group $\CH(X)\otimes_{\mathbb{Z}} \mathbb{F}$ of cycles on $X$ modulo rational equivalence with coefficients in a finite field $\mathbb{F}$. We write $\Ch(\overline{X})$ for the colimit of all $\Ch(X_K)$, where $K$ runs through all field extensions $K/F$ and if $X$ is integral we denote by $F(X)$ its function field.

For any field extension $L/F$, an element lying in the image of the natural morphism of $\Ch(X_{L})\longrightarrow \Ch(\overline{X})$ is called $L$-rational. The image of any correspondence $\alpha\in \Ch(X_{L})$ under the canonical morphism $\Ch(X_{L})\longrightarrow\Ch(\overline{X})$ is denoted by $\overline{\alpha}$. An $F$-variety $X$ is geometrically split if the Grothendieck Chow motive of $X_{\overline{F}}$=$X\times_{\Spec(F)}\Spec(\overline{F})$ with coefficients in $\mathbb{F}$ is isomorphic to a finite direct sum of Tate motives, for an algebraic closure $\overline{F}/F$. The variety $X$ satisfies the Rost nilpotence principle with coefficients in $\mathbb{F}$ if for any field extensions $L/E/F$ the kernel of the restriction map $res_{L/E}:\End(M(X_E))\longrightarrow \End(M(X_L))$ consists of nilpotents.

As shown in \cite{chermergil}, any projective homogeneous $F$-variety under the action of a semisimple affine algebraic group is geometrically split and satisfies the Rost nilpotence principle. It follows by \cite[Corollary 35]{chermer} (see also \cite[Corollary 2.6]{upper}) that the Grothendieck Chow motive of these varieties with coefficients in $\mathbb{F}$ decomposes in an essentially unique way as a direct sum of indecomposable motives. The study of these decompositions have already shown to be very fruitful (see \cite{hoi}, \cite{upper}, \cite{vish}).

The notion of \emph{upper} motives, previously defined by Vishik in the context of quadrics in \cite{vish}, was further developed by Karpenko in \cite{upper} to describe the indecomposable motives lying in the motivic decomposition of projective homogeneous varieties. If $X$ is a homogeneous $F$-variety, $E/F$ a field extension and the upper motive of $M(X_E)$ is a direct summand of another motive $M_E$, \cite[Theorem 4.15]{vish} and \cite[Proposition 4.6]{hoi} give sufficient conditions for the upper motive of $X$ to be a direct summand of $M$. The purpose of the present note is to push these ideas further. We define the notions of \emph{upper}, \emph{lower} and \emph{outer} direct summands of a direct summand $N$ of the motive of a geometrically split $F$-variety. We then show some lifting property of outer summands of $N_E$ to the base field with the following result.

\begin{Thm}\label{main}
Let $N$ be a direct summand of the motive (with coefficients in $\mathbb{F}$) of a geometrically split, geometrically integral $F$-variety $X$ satisfying the Rost nilpotence principle with coefficients in $\mathbb{F}$ and $M$ a twisted direct summand of the motive of another $F$-variety $Y$. Assume that there is a field extension $E/F$ such that
\begin{enumerate}
\item every $E(X)$-rational cycle in $\Ch(\overline{X\times Y})$ is $F(X)$-rational;
\item the motive $N_E$ has an indecomposable outer direct summand which is also a direct summand of the motive $M_E$.
\end{enumerate}
Then the motive $N$ has an outer direct summand which is also a direct summand of $M$.
\end{Thm}

Theorem \ref{main} allows one to descend outer motives of direct summands projective homogeneous varieties which appear on some field extension $E/F$ of the base field. This generalizes \cite[Proposition 4.6]{hoi}, one of the key ingredients in the proof of \cite[Theorem 1.1]{hoi}, replacing the whole motive of a variety $X$ by a direct summand. To replace $X$ by an arbitrary direct summand, one needs to construct explicitly the rational cycles to get an outer summand defined over $F$, and thus theorem \ref{main} gives a new proof of \cite[Proposition 4.6]{hoi}. Note that assumption 1 of theorem \ref{main} holds if the field extension $E(X)/F(X)$ is unirational, i.e. if there is a field extension $L/E(X)$ such that $L/F(X)$ is purely transcendental.

The following particular case of theorem \ref{main} was used by Garibaldi, Petrov and Semenov in \cite{shells} to both determine all the motivic decompositions of homogeneous $F$-varieties of inner type $E_6$ and prove a conjecture of Rost and Springer in \cite{shells}. 

\begin{coroll}(\cite[{Proposition 3.2}]{shells})\label{coroll} Let $X$ and $Y$ be two projective homogeneous $F$-varieties for a semisimple affine algebraic group, and let $M$ and $N$ be direct summands of the motives of $Y$ and $X$ respectively with coefficients in $\mathbb{F}$. Assume that $N_{F(Y)}$ is an indecomposable direct summand of $M_{F(Y)}$ and $Y$ has an $F(X)$-point. Then $N$ is a direct summand of $M$.
\end{coroll}
\begin{proof}[Proof]Setting $E=F(Y)$, the field extension $E(X)/F(X)$ is purely transcendental, hence assumption $1)$ of theorem \ref{main} holds.
\end{proof}

\section{Grothendieck Chow motives}

Our main reference for the construction of the category of Grothendieck Chow motives over $F$ with coefficients in $\mathbb{F}$ is \cite[$\S$63-65]{EKM}.

Let $X$ and $Y$ be two $F$-varieties and $X=\coprod_{k=1}^nX_k$ be the decomposition of $X$ as a disjoint union of irreducible components with respective dimension $d_1$,...,$d_n$. For any integer $i$ the group of \emph{correspondences} between $X$ and $Y$ of degree $i$ with coefficients in $\mathbb{F}$ is defined by $\Corr_i(X,Y)=\coprod_{k=1}^n\Ch_{d_k+i}(X_k\times Y).$
We now consider the category $\CCC(F;\mathbb{F})$ whose objects are pairs $X[i]$, where $X$ is an $F$-variety and $i$ is an integer. Morphisms are defined in terms of correspondences by $\Hom_{\CCC(F;\mathbb{F})}(X[i],Y[j])=\Corr_{i-j}(X,Y)$. For any correspondences $f:X[i]\rightsquigarrow Y[j]$ and $g:Y[j]\rightsquigarrow Z[k]$ in $\Mor(\CCC(F;\mathbb{F}))$ the composite $g\circ f:X[i]\rightsquigarrow Z[k]$ is defined by $$g\circ f=\left(^X\!p^Z_Y\right)_{\ast}\left((^{X\times Y}\!\!p_Z)^{\ast}(f)\cdot (p^{Y\times Z}_X)^{\ast}(g)\right)~~~~~~~~~\mbox{($\ast$)}$$
where $^U\!p^W_V:U\times V\times W\rightarrow U\times W$ is the natural projection.

The category $\CCC(F;\mathbb{F})$ is preadditive and its additive completion $\CR(F;\mathbb{F})$ is the category of correspondences over $F$ with coefficients in $\mathbb{F}$, which has a structure of tensor additive category given by $X[i]\otimes Y[j]=(X\times Y)[i+j]$.
The category $\CM(F;\mathbb{F})$ of Grothendieck Chow motives with coefficients in $\mathbb{F}$ is the pseudo-abelian envelope of the category $\CR(F;\mathbb{F})$. Its objects are couples $(X,\pi)$, where $X$ is an object of the category $\CR(F;\Lambda)$, and $\pi\in \End(X)$ is a projector (i.e. $\pi \circ \pi=\pi$). Morphisms are given by $\Hom_{\CM(F;\mathbb{F})}\left((X,\pi),(Y,\rho)\right)=\rho\circ \Hom_{\CR(F;\mathbb{F})}(X,Y)\circ \pi$ and the objects of $\CM(F;\mathbb{F})$ are called \emph{motives}. For any $F$-variety $X$ the motives $(X[i],\Gamma_{id_X})$ (where $\Gamma_{id_X}$ is the graph of the identity of $X$) will be denoted $X[i]$ and $X[0]$ is the motive of $X$. The motives $\mathbb{F}[i]=\Spec(F)[i]$ are the \emph{Tate motives}. 

\begin{lemma}\label{sousproj}Let $(X,\pi)$ be a direct summand of the motive of an $F$-variety $X$. A motive $M$ is a direct summand of $(X,\pi)$ if and only if $M$ is isomorphic to $(X,\rho)$, for some projector $\rho$ satisfying $\pi\circ \rho \circ \pi=\rho$.
\end{lemma}
\begin{proof}Since $\End\left((X,\pi)\right)=\pi\circ \Ch_{\dim(X)}(X\times X)\circ \pi$, any projector $\rho$ in $\End\left((X,\pi)\right)$ satisfies $\pi\circ \rho\circ \pi=\rho$.
\end{proof}

\begin{Df}Let $M\in \CM(F;\mathbb{F})$ be a motive and $i$ an integer. The $i$-dimensional Chow group $\Ch_i(M)$ of $M$ is defined by $\Hom_{\CM(F;\mathbb{F})}(\mathbb{F}[i],M)$. The $i$-codimensional Chow group $\Ch^i(M)$ of $M$ is defined by $\Hom_{\CM(F;\mathbb{F})}(M,\mathbb{F}[i])$.
\end{Df}
For any field extension $E/F$ and any correspondence $\alpha:X[i]\rightsquigarrow Y[j]$ the pull-back of $\alpha$ along the natural morphism $(X\times Y)_E\rightarrow X\times Y$ will be denoted $\alpha_E$. If $N=(X,\pi)[i]$ is a twisted motivic direct summand of $X$, the motive $(X_E,\pi_E)[i]$ will be denoted $N_E$.

Finally the category $\CM(F;\mathbb{F})$ is endowed with a duality functor. If $X$ and $Y$ are two $F$-varieties and $\alpha\in \Ch(X\times Y)$ is a correspondence, the image of $\alpha$ under the exchange isomorphism $X\times Y \rightarrow Y\times X$ is denoted $^t\alpha$. The \emph{duality functor} is the additive functor $\dag:\CM(F;\Lambda)^{op}\longrightarrow \CM(F;\Lambda)$ determined by the formula $M(X)[i]^{\dag}=M(X)[-\dim(X)-i]$ and such that for any correspondence $\alpha:X[i]\rightsquigarrow Y[j]$, $\alpha^{\dag}=^t\!\alpha$.

\section{Direct summands of geometrically split $F$-varieties}\label{deftop}

Throughout this section we consider a geometrically split $F$-variety $X$ and $E/F$ a splitting field of $X$. By \cite[Proposition 1.5]{requivtori} the pairing
$$\fonction{\Psi}{\Ch(X_E)\times \Ch(X_E)}{\mathbb{F}}{(\alpha,\beta)}{deg(\alpha\cdot\beta)}$$
is non degenerate hence gives rise to an isomorphism of $\mathbb{F}$-modules between $\Ch(X_E)$ and its dual space $\Hom_{\mathbb{F}}(\Ch(X_E),\mathbb{F})$ given by $\alpha\mapsto \Psi(\alpha,\cdot)$. The dual basis of a homogeneous basis $(x_k)_{k=1}^n$ of $\Ch(X_E)$ with respect of $\Psi$ is the basis $(x_k^{\ast})_{k=1}^n$ of $\Ch(X_E)$ such that for any $1\leq i,j\leq n$, $\Psi(x_i,x^{\ast}_j)=\delta_{ij}$, where $\delta_{ij}$ is the Kronecker symbol. By definition of the composition ($\ast$) in $\CM(F;\mathbb{F})$, if $y$ (resp. $y'$) lies in $\Ch(Y)$ (resp. $\Ch(Y')$) for two other $F$-varieties $Y$ and $Y'$ and if $(i,j)$ are two integers, the composition of the correspondences $x_i\times y\in \Ch(X_E\times Y)$ and $y'\times x^{\ast}_j\in \Ch(Y'\times X_E)$ is given by
\begin{align}\label{formula}(x_i\times y)\circ (y'\times x_j^{\ast})=\delta_{ij}(y'\times y)\in \Ch(Y'\times Y).\end{align}
Note that the Kunneth decomposition holds in $\Ch(X_E\times Y)$ and $\Ch(Y'\times X_E)$ in the view of \cite[Proposition 64.3]{EKM}, since $X_E$ is split, and thus the cycles of $\Ch(X_E\times Y)$ and $\Ch(Y'\times X_E)$ may always be written that way.\\

\textbf{Upper, lower and outer motives.} Let $\pi\in \Ch_{\dim(X)}(X\times X)$ be a non-zero projector and $N$$=$$(X,\pi)$ the associated summand of the motive of $X$. The \emph{base} of $N$ is the set $\mathcal{B}(N)=\{i\in \mathbb{Z},~\Ch_i(N_E)\mbox{ is not trivial}\}$. The \emph{bottom} of $N$ (denoted by $b(N)$) is the least integer of $\mathcal{B}(N)$ and the \emph{top} of $N$ (denoted by $t(N)$) is the greatest integer of $\mathcal{B}(N)$. We now introduce the notion of upper and lower direct summands of $N$, previously introduced by Vishik in the context of the motives of quadrics in \cite[Definition 4.6]{vish}.
 
\begin{Df}\label{factext}Let $N$ be a direct summand of the twisted motive of a geometrically split $F$-variety and $M$ a motivic direct summand of $N$. We say that
\begin{enumerate}
\item $M$ is \emph{upper} in $N$ if $b(M)=b(N)$;
\item $M$ is \emph{lower} in $N$ if $t(M)=t(N)$;
\item $M$ is \emph{outer} in $N$ if $M$ is both lower and upper in $N$. 
\end{enumerate}
\end{Df}

\begin{Rmk}\label{rmk}Keeping the same $F$-variety $X$ and any direct summand $N=(X,\pi)$, consider a homogeneous basis $(x_k)_{k=1}^n$ of $\Ch(X_E)$ and its dual basis $(x_k^{\ast})_{k=1}^n$. The base, bottom and top of $N$ can be easily determined by the decomposition 
$$\pi_E=\sum_{i,j=1}^n\pi_{i,j}(x_i\times x_j^{\ast})$$
noticing that $\mathcal{B}(N)=\{\dim(x_i),~\mbox{$\pi_{i,j}\neq 0$ for some $j$}\}$. 
\end{Rmk}

\begin{lemma}\label{dag}Let $N$ be a motivic direct summand of a geometrically split $F$-variety and $M$ a direct summand of $N$. Then $M$ is lower in $N$ (resp. upper in $N$) if and only if the dual motive $M^{\dag}$ is upper in $N^{\dag}$ (resp. $M^{\dag}$ is lower in $N^{\dag}$).
\end{lemma}
\begin{proof}[Proof]For any motive $O$ and for any integer $i$, $\Ch^i(O^{\dag})=\Ch_{-i}(O)$. It follows that $b(O^{\dag})=-t(O)$ and $t(O^{\dag})=-b(O)$.
\end{proof}

\textbf{The Krull-Schmidt property.} Let $\CCC$ be a pseudo-abelian category and $\mathfrak{C}$ be the set of the isomorphism classes of objects of $\CCC$. We say that the category $\CCC$ satisfies the \emph{Krull-Schmidt property} if the monoid $(\mathfrak{C},\oplus)$ is free. The Krull-Schmidt property holds for the motives of geometrically split $F$-varieties satisfying the Rost nilpotence principle in $\CM(F;\mathbb{F})$ by \cite[Corollary 2.6]{upper}.\\

\textbf{Proof of the main result.} In order to prove theorem \ref{main}, we will need the following lemma, which will allow us to construct explicitly the rational cycles lifting outer motives to the base field.

\begin{lemma}\label{lmme}Let $N$ be a motivic direct summand of a geometrically split, geometrically irreducible $F$-variety $X$ satisfying the Rost nilpotence principle and $M$ a twisted direct summand of an $F$-variety $Y$. Assume the existence of a field extension $E/F$ such that
\begin{enumerate}
\item any $E(X)$-rational cycle in $\Ch(\overline{X\times Y})$ is $F(X)$-rational;
\item there are two correspondences $\alpha:N_E\rightsquigarrow M_E$ and $\beta:M_E\rightsquigarrow N_E$ such that $\beta\circ \alpha$ is a projector and $(X_E,\beta\circ \alpha)$ is a lower direct summand of $N_E$.
\end{enumerate}
Then there are two correspondences $\gamma:N\rightsquigarrow M$ and $\delta:M_E\rightsquigarrow N_E$ such that $(X_E,\delta\circ \gamma_E)$ is a direct summand of $N_E$ which contains all lower indecomposable direct summands of $(X_E,\beta\circ \alpha)$. Furthermore if $\overline{\beta}$ is $F$-rational, then $\overline{\delta}$ is also $F$-rational.
\end{lemma}
\begin{proof}We may assume by lemma \ref{sousproj} that $M$=$(Y,\rho)[i]$ and $N$=$(X,\pi)$. We construct explicitly the two correspondences $\gamma$ and $\delta$. Since $E(X)$ is a field extension of $E$, $\overline{\alpha}$ is $E(X)$-rational, hence $F(X)$-rational by assumption $1$. Let $L/F$ be a field extension. Then the morphism $\Spec(F(X_L))\longrightarrow X_L$ induces a pull-back morphism $\varepsilon^{\ast}:\Ch(\overline{X\times Y\times X})\longrightarrow \Ch(\overline{(X \times Y)_{F(X)}})$ which maps $F$-rational cycles onto $F(X)$-rational cycles by \cite[Corollary 57.11]{EKM}, and so there is a cycle $\alpha_1\in \Ch(X\times Y\times X)$ such that $\varepsilon^{\ast}(\overline{\alpha_1})=\overline{\alpha}$. Since $\varepsilon^{\ast}$ maps any homogeneous cycle $\sum_ix_i\times y_i\times 1$ to $\sum_i x_i\times y_i$ and vanishes on homogeneous cycles whose codimension on the third factor is strictly positive, we have $\overline{\alpha_1}=\overline{\alpha}\times 1+\cdots$ where $"\cdots"$ is a linear combination of homogeneous cycles in $\Ch(\overline{X\times Y\times X})$ with strictly positive codimension on the third factor.

We now look at $\alpha_1$ as a correspondence $X\rightsquigarrow Y\times X$ and consider the cycle $\alpha_2=\alpha_1\circ\pi$. By formula \ref{formula} we have $$\overline{\alpha_2}=\left(\overline{\alpha}\times 1\right)\circ \overline{\pi}+\cdots,$$ where $"\cdots"$ is a linear combination of homogeneous cycles in $\Ch(\overline{X \times Y\times X})$ with dimension at most $t(N)$ on the first factor (since these terms come from the first factors of $\overline{\pi}$) and strictly positive codimension on the third factor (since these terms come from the third factors of $\overline{\alpha_1}-\overline{\alpha}\times 1$). Finally considering the pull-back of the morphism $\Delta:X\times Y\rightarrow  X\times Y\times X$ induced by the diagonal embedding $X$ and setting $\alpha_3=\Delta^{\ast}(\alpha_2)$, we have
$$\overline{\alpha_3}=\overline{\alpha}\circ\overline{\pi}+\cdots$$
where $"\cdots"$ stands for a linear combination of homogeneous cycles in $\Ch(\overline{X \times Y})$ with dimension strictly lesser than $t(N)$ on the first factor.

Composing with $\overline{\pi}\circ \overline{\beta}$ on the left and $\overline{\pi}$ on the right, we get that
$$\overline{\pi}\circ \overline{\beta} \circ\overline{\alpha_3}\circ\overline{\pi}=\overline{\pi}\circ \overline{\beta}\circ\overline{\alpha}\circ \overline{\pi}+\xi$$
where $\xi$ is a linear combination of homogeneous cycles of strictly lesser dimension than $t(N)$ on the first factor since they come from the first factors of $\overline{\alpha_3}\circ \overline{\pi}- \overline{\alpha}\circ \overline{\pi}$. The correspondence $\beta\circ \alpha$ defines a direct summand of the motive $N_E$ and thus by lemma \ref{sousproj}
$$\overline{\pi}\circ \overline{\beta} \circ\overline{\alpha_3}\circ\overline{\pi}=\overline{\beta}\circ\overline{\alpha}+\xi.$$

By formula \ref{formula}, $\overline{\beta}\circ \overline{\alpha} \circ \xi$, $\xi\circ \xi$ and $\xi\circ \overline{\beta}\circ \overline{\alpha}$ are linear combinations of homogeneous cycles of dimension strictly lesser than $t(N)$ on the first factor. Repeating the same procedure and since $k\circ h$ is a projector, we see that for any integer $n$

\begin{align}\label{finn}\left(\overline{\pi}\circ\overline{\beta}\circ \overline{\alpha_3}\circ\overline{\pi}\right)^{n}=\overline{\beta}\circ \overline{\alpha}+\cdots
\end{align}
where $"\cdots"$ is a linear combination of homogeneous cycles in $\Ch(\overline{X\times X})$ with dimension on the first factor strictly lesser than $t(N)$. Since the direct summand $(X,\beta\circ \alpha)$ is lower, all these correspondences are non-zero and by \cite[Corollary 2.2]{upper} an appropriate power $\left(\pi_E\circ \beta\circ (\alpha_3)_E\circ\pi_E\right)^{\circ n_0}$ is a projector. If we set $\gamma$$=$$\rho\circ \alpha_3\circ \pi$ and $\delta$$=$$\left(\pi_E\circ \beta\circ (\alpha_3)_E\circ\pi_E\right)^{\circ n_0-1}\circ \pi_E\circ \beta$, we see that $\overline{\delta}$ is $F$-rational if $\overline{\beta}$ is $F$-rational. The correspondence $\delta\circ \gamma_E$ is a projector which defines a direct summand of $N_E$ by lemma \ref{sousproj}.

Consider the decomposition $\overline{\beta}\circ \overline{\alpha}$$=$$\sum_{i,j=1}^sp_{ij}(x_i\times x_j^{\ast})$ of $\overline{\beta}\circ\overline{\alpha}$ with respect to a basis $(x_i)_{i=1}^s$ of $\Ch(\overline{X})$. By formula \ref{finn}, the decomposition of $\overline{\delta}\circ \overline{\gamma}$ in $(x_i\times x_j^{\ast})_{i,j=1}^s$ has a non-zero coefficient for any couple $(i,j)$ such that $p_{ij}$ is non zero and $\dim(x_i)=t((X_E,\beta\circ \alpha))$. The Krull-Schmidt property and remark \ref{rmk} then imply that any lower indecomposable direct summand of $(X_E,\beta\circ \alpha)$ is a direct summand of $(X,\delta\circ \gamma_E)$.
\end{proof}

We now show how we can derive the proof of theorem \ref{main} from the rational cycles constructed in lemma \ref{lmme}. To lift the outer motive to the base field $F$, we apply lemma \ref{lmme} and the duality functor twice in order produce two correspondences which are defined on the base field.

\begin{proof}[Proof of Theorem 1.]Let $O$=$(X_E,\kappa)$ be an outer indecomposable direct summand of $N_E$ which is also a direct summand of $M_E$. We prove theorem \ref{main} by applying lemma \ref{lmme} once, then the duality functor and finally lemma \ref{lmme} another time to get all our correspondences defined over the base field $F$.

Since $O$ is a direct summand of $M_E$, there are two correspondences $\alpha:N_E\rightsquigarrow M_E$ and $\beta:M_E\rightsquigarrow N_E$ such that $\beta\circ \alpha=\kappa$. Moreover $O$ is lower in $N_E$, so lemma \ref{lmme} justifies the existence of two other correspondences $\alpha ':N\rightsquigarrow M$ and $\beta ':M_E\rightsquigarrow N_E$ such that $O_2=(X_E,\beta '\circ \alpha '_E)$ is a direct summand of $N_E$, and the motive $O_2$ is outer in $N_E$ since it contains $O$. The dual motive $O_2^{\dag}=(X_E,^t\!\alpha '_E\circ ^t\!\beta ')[-\dim(X)]$ is therefore outer in $N_E^{\dag}$ by lemma \ref{dag} and is a direct summand of the dual motive $M_E^{\dag}$. Twisting these three motives by $\dim(X)$, we can apply lemma \ref{lmme} again. The correspondence $\overline{^t\alpha '}$ is $F$-rational, so lemma \ref{lmme} gives two correspondences $\gamma:N^{\dag}\rightsquigarrow M^{\dag}$ and $\delta:M^{\dag}\rightsquigarrow N^{\dag}$ such that the motive $(X_E,\delta_E\circ \gamma_E)$ is both an outer direct summand of $N^{\dag}$ (since it contains the dual motive $O^{\dag}$) and a direct summand of $M^{\dag}$. Transposing again, the motive $(X,^t\!\gamma\circ ^t\!\delta)$ is an outer direct summand of $N$ and a direct summand of $M$.
\end{proof}

\section{Motivic decompositions for groups of inner type $E_6$}

The purpose of this section is to discuss the complete classification of the motivic decompositions of projective homogeneous varieties of inner type $E_6$, which is achieved in \cite{shells}. Let $G$ be an algebraic group of inner type $E_6$ and $X$ a projective $G$-homogeneous variety. We choose the following numbering of the Dynkin diagram $G$.\\
\vspace{0,5cm}

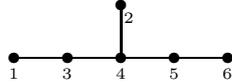
\makebox[16cm][c]{
\begin{picture}(100,7)
\put(40,20){\circle*{4}}
\multiput(0,0)(20,0){5}{\circle*{4}}
\put(-2,-8){$\scriptscriptstyle 1 $}
\put(18,-8){$\scriptscriptstyle 3 $}
\put(38,-8){$\scriptscriptstyle 4 $}
\put(58,-8){$\scriptscriptstyle 5 $}
\put(78,-8){$\scriptscriptstyle 6 $}
\put(41,13){$\scriptscriptstyle 2 $}

\put(40,0){\line(0,1){20}}
\put(0,0){\line(1,0){80}}

\end{picture}
}

\vspace{0,5cm}

The results of \cite{Jinv} show that in the case where $X$ is \emph{generically split}, any indecomposable summand of the $\mathbb{F}_p$-motive of $X$ is isomorphic to a shift of the upper motive $\mathcal{R}_p(G)$ of the variety of Borel subgroups of $G$. Furthermore, the structure of the motives $\mathcal{R}_p(G)$ is determined in \cite{Jinv} in terms of the so-called $J$-invariant modulo $p$ of $G$.

The $J$-invariant was first introduced by Vishik in \cite{vish4} in the context of quadratic forms. Petrov, Semenov and Zainoulline define in \cite{Jinv} the notion of $J$-invariant modulo $p$ of an arbitrary semisimple algebraic group $G$, denoted by $J_p(G)$, which is an $r$-tuple of integers $(j_1,...,j_r)$ given by the rational cycles in $\Ch(\overline{G})$. By \cite[Table 4.13]{Jinv}, the $J$-invariant modulo $3$ of a semisimple adjoint algebraic group of inner type $E_6$ is $(j_1,j_2)$, with $0\leq j_1\leq 2$ and $0\leq j_2\leq 1$.

Another invariant attached to $G$ is the Tits index, which consists of the data of the Dynkin diagram of $G$ with some vertices being circled. The complete classification of the Tits indices of type $E_6$, provided in \cite{tits}, is as follows :\\
\vspace{0,5cm}

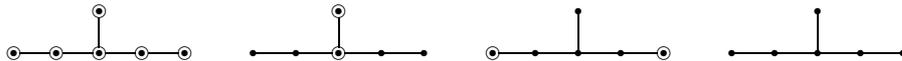
\makebox[16cm][c]{
\scalebox{0.8}{
\begin{tabular}{cccc}

\begin{picture}(100,7)
\put(40,20){\circle*{3}}
\put(40,20){\circle{6}}
\multiput(0,0)(20,0){5}{\circle*{3}}
\multiput(0,0)(20,0){5}{\circle{6}}

\put(40,3){\line(0,1){14}}
\put(3,0){\line(1,0){14}}
\put(23,0){\line(1,0){14}}
\put(43,0){\line(1,0){14}}
\put(63,0){\line(1,0){14}}
\end{picture}
&
\begin{picture}(100,7)
\put(40,20){\circle*{3}}
\put(40,20){\circle{6}}
\multiput(0,0)(20,0){5}{\circle*{3}}
\put(40,0){\circle{6}}

\put(40,3){\line(0,1){14}}
\put(0,0){\line(1,0){20}}
\put(20,0){\line(1,0){17}}
\put(43,0){\line(1,0){17}}
\put(60,0){\line(1,0){20}}
\end{picture}
&
\begin{picture}(100,7)
\put(40,20){\circle*{3}}
\multiput(0,0)(20,0){5}{\circle*{3}}
\multiput(0,0)(80,0){2}{\circle{6}}

\put(40,0){\line(0,1){20}}
\put(3,0){\line(1,0){17}}
\put(20,0){\line(1,0){20}}
\put(40,0){\line(1,0){20}}
\put(60,0){\line(1,0){17}}
\end{picture}
&
\begin{picture}(100,7)
\put(40,20){\circle*{3}}
\multiput(0,0)(20,0){5}{\circle*{3}}

\put(40,0){\line(0,1){20}}
\put(0,0){\line(1,0){20}}
\put(20,0){\line(1,0){20}}
\put(40,0){\line(1,0){20}}
\put(60,0){\line(1,0){20}}
\end{picture}
\end{tabular}
}
}

\vspace{0,5cm}

Let $\Theta$ be a subset of the vertices of the Dynkin diagram of $G$ and $X_{\Theta}$ a projective $G$-homogeneous variety of type $\Theta$. By \cite[Table 3.6]{Jinv} and a case by case anylisis of the above Tits indices, the variety $X_{\Theta}$ is either split or generically split if $p\neq 3$, $\Theta\neq \{2\}$, $\{4\}$ or $\{2,4\}$ and if $j_1=0$. Finally, the results of \cite{chermergil} and \cite[\S 8]{shells} imply that the understanding of the motivic decompositions of $X_{\{2\}}$, $X_{\{4\}}$, and $X_{\{2,4\}}$ is reduced to the study of the upper motive of $X_{2}$ in $\CM(F;\mathbb{F}_3)$, which is denoted by $M_{j_1,j_2}$.

Using theorem \ref{main}, Garibaldi, Petrov and Semenov provide some restrictions on the Poincaré polynomial of the upper motive of $X$, if $X$ is an anisotropic projective homogeneous variety satisfying some technical assumptions (see \cite[Proposition 7.6]{shells}). Assuming that $J_3(G)$$=$$(1,0)$, they observe that although the variety $X_{\{2\}}$ satisfies all those technical assumptions, the Poincaré polynomial of its upper motive does not match with the conclusion of \cite[Proposition 7.6]{shells}. In particular $X_{\{2\}}$ has a $0$-cycle of degree coprime to $3$, and $M_{1,0}$ is the Tate motive $\mathbb{F}_3$.

Furthermore, the authors deduce from the fact that $M_{1,0}$ is the Tate motive that the $J$-invariant modulo $3$ of $G$ cannot be $(2,0)$ (see \cite[Corollary 8.10]{shells}). Indeed, if $J_3(G)$$=$$(2,0)$ and $\SB(3,A)$ is the Severi-Brauer variety of right ideals of reduced dimension $3$ in the Tits algebra of $G$, then $J_3(G_{F(\SB(3,A))})$$=$$(1,0)$. In particular, $X_2$ has a zero-cycle of degree coprime to $3$ over the function field of $\SB(3,A)$. It follows that the upper motive of $X_2$ would be isomorphic to the upper motive of $\SB(3,A)$, and thus the canonical $3$-dimensions of $X_2$ and $\SB(3,A)$ would be equal, a contradiction.\\

The authors use similar techniques to provide isotropy criteria for projective homogeneous varieties. They consider several varieties which satisfy the technical assumptions of \cite[Proposition 7.6]{shells} without fulfilling its conclusion, and thus have a zero cycle of degree $3$ (or a rational point). These examples include projective homogeneous varieties for orthogonal group with application to the isotropy of varieties of type $E_7$ and varieties of type $E_8$ (see \cite[Lemma 10.15,10.21]{shells}). They also produce with these techniques isotropy criteria for projective homogeneous varieties in terms of the Rost invariant (see \cite[Proposition 10.18]{shells} for type $E_7$ and \cite[Propositions 10.22]{shells} for type $E_8$).\\

\textbf{Acknowledgements :} I am grateful to N. Karpenko for raising this question and for his suggestions. I also would like to thank N. Semenov for the very useful conversations.

\begin{footnotesize}

\end{footnotesize}

\end{document}